\documentclass[runningheads]{llncs}
\usepackage[utf8]{inputenc}
\usepackage{graphics}
\usepackage{mathtools}
\usepackage{amsfonts}
\usepackage{multirow}
\usepackage{multicol}
\usepackage{enumerate}
\usepackage{MnSymbol}
\usepackage{graphicx}
\usepackage[overload]{empheq}
\usepackage{bbm}
\usepackage{hyperref}

\usepackage[english]{babel}

\usepackage{tkz-berge}
\definecolor{iceberg}{rgb}{0.44, 0.65, 0.82}

\usepackage{csquotes}

\newtheorem{PR}[theorem]{Proposition}

\begin{document}

\author{Olha Silina}
\title{Covering the edges of a graph with perfect matchings}
\institute{Carnegie Mellon University, Pittsburgh, PA 15213, USA \\
\email{osilina@andrew.cmu.edu}\\
}

\maketitle

%\author{Olha Silina\footnote{The author is a graduate student at the Department of Mathematical Sciences at Carnegie Mellon University. Email: \texttt{osilina@andrew.cmu.edu}}}
%\date{\today}
%\maketitle

\begin{abstract}
    
 An $r$-graph is an $r$-regular graph with no odd cut of size less than $r$. A well-celebrated result due to Lov\'asz says that for such graphs the linear system $Ax = \mathbbm{1}$ has a solution in $\mathbb{Z}/2$, where $A$ is the $0,1$ edge to perfect matching incidence matrix. Note that we allow $x$ to have negative entries. In this paper, we present an improved version of Lov\'asz's result, proving that, in fact, there is a solution $x$ with all entries being either integer or $+1/2$ and corresponding to a linearly independent set of perfect matchings. Moreover, the total number of $+1/2$'s is at most $6k$, where k is the number of Petersen bricks in the tight cut decomposition of the graph.

%\keywords{Perfect matching  \and Lattice .}
 \end{abstract}

\section{Introduction}

Let $G=(V,E)$ be a connected nontrivial graph. A set of edges $M \subseteq E$ is a \emph{perfect matching} if every vertex is incident with exactly one edge of $M$. A graph is \emph{matching-covered} if every edge belongs to a perfect matching. Matching theory is a rich and widely studied area, with extensive work done around matching-covered graphs, starting with \cite{Edmonds1965MaximumMA,Edmonds1982BrickDA,lovaszplummer}. A particular class of matching-covered graphs are so-called \emph{$r$-graphs}: a graph is an \emph{$r$-graph} if every vertex has degree $r$ and all odd cuts have size at least $r$ (a \emph{cut} is the set of edges with exactly one endpoint in $S$ for some set $S \subseteq V$; a cut is \emph{odd} if $|S|$ is odd).

The notion of $r$-graphs was first introduced by Seymour, who also showed that an $r$-graph is always matching-covered \cite{Seymour}. Being an $r$-graph is a necessary condition for the graphs whose edge sets can be partitioned into $r$ perfect matchings. One question that naturally arises here is whether this is a sufficient condition, i.e. whether an $r$-graph contains $r$ disjoint perfect matchings. Unfortunately, it is already false for $r=3$ in the case of Petersen graph, which does not contain two disjoint perfect matchings. However, there are several interesting relaxations of this problem. One of the central results in this area is due to Lov\'asz, who got a characterization for the lattice $\mathcal{L} = \{\sum \alpha_i A_i : \alpha_i \in \mathbb{Z}\}$ generated by the incidence vectors of perfect matchings $A_i$ of a matching-covered graph $G$ using dual lattice theory \cite{lovaszPM}.
 It follows from the characterization above that for an $r$-graph, the vector of all-twos $2\cdot\mathbbm{1}$ belongs to $\mathcal{L}$. Hence, the vector of all-ones $\mathbbm{1}$ can be obtained as a half-integral (but not necessarily non-negative) combination of the perfect matchings of $G$. In this paper, we strengthen this result by proving the existence of a half-integral solution in which there are only a few fractional coefficients, all equal to exactly $+1/2$.

More formally, given an $r$-graph $G = (V,E)$ with $n$ vertices and $m$ edges (notice $2m = nr$), we want to study the set of solutions to $Ax=\mathbbm{1}$, where $A$ is a $\{0,1\}$ incidence matrix with rows corresponding to the edges of $G$ and columns corresponding to the perfect matchings in $G$ (we will adopt this notation throughout the paper). We consider graphs with no loops, but parallel edges are allowed.

The main result of this paper is:
\begin{theorem}[Main theorem]\label{main}
Let $G$ be an $r$-graph with $m$ edges and $n$ vertices and let $A$ be its edge to perfect matching incidence matrix. Then there is a solution $x^*$ to $Ax = \mathbbm{1}$ satisfying the following conditions: \begin{enumerate}
    \item all non-integral entries of $x^*$ are equal to $+1/2$;
    \item $x^*$ has at most $m-n+1$ non-zero entries;
    \item all non-zero entries of $x^*$  correspond to a linearly independent set of perfect matchings of $G$;
    \item the total number of $+1/2$'s in $x^*$ is at most $6p$, where $p$ is the number of Petersen bricks of $G$.
\end{enumerate}
\end{theorem}

The notions of a brick and the Petersen brick are defined later, but for a reference, $p$ is bounded from above by the total number of vertices, $n$. We use $+1/2$ instead of $1/2$ to highlight the fact that these fractional entries are positive. Intuitively, there are two obstacles to a solution to $Ax = \mathbbm{1}$ with $x$ being a $\{0,1\}$ vector: one is that $x$ has to be fractional, and the other one is that $x$ has to have negative entries. Our result guarantees a solution with all entries violate at most one of the above conditions.
%Notice that the matrix $A$ can have exponentially many columns as a function of the input (the number of edges and vertices in the graph), so this problem cannot be solved efficiently using standard linear algebra methods.

The rest of this paper is organized as follows. In Section \ref{sect:tight_dicuts} we define a central to the whole argument definition of \textit{tight} cuts and explain a process of deconstructing a perfect-matching covered graph into basic pieces, called \textit{bricks} and \textit{braces}. Section \ref{sect:general_basis_combination} summarizes a classical argument (see, e.g. \cite{Murty}) for combining two solutions of a smaller problem into a solution to the original problem. An improved version of this process together with an inductive analysis for the support size and the total number of fractional entries in the solution are presented in Section \ref{sect:main}. It will then suffice to prove the Theorem \ref{main} for the case of bricks and braces. The latter case is fairly straightforward, while the former case will be addressed in Sections \ref{sect:petersen} and \ref{sect:non-petersen} for the cases of Petersen and non-Petersen bricks, respectively.

\section{Tight cut decomposition}\label{sect:tight_dicuts}

Tight cut decomposition plays an essential role in many results regarding perfect matchings. We begin by introducing a special class of odd cuts. Let $G=(V,E)$ be a perfect-matching covered graph.

\begin{definition}
    A \emph{tight cut} is an odd cut $C$($\subseteq E$) such that every perfect matching intersects $C$ in exactly one edge.
A cut is \emph{trivial} if one of its shores is a single vertex.
Otherwise, a cut is \emph{non-trivial}.
\end{definition}

It is clear that all trivial cuts of $G$ are tight since every perfect matching has exactly one edge incident to every vertex. However, in general there could be other non-trivial tight cuts. Let us prove the following statement:
\begin{lemma}\label{lemma1}
In a matching-covered graph $G$ with $Ax=\mathbbm{1}$ feasible, all tight cuts have the same size. In particular, $G$ is regular.
\end{lemma}
\begin{proof}
    Consider any tight cut $C$ and its characteristic vector $\mathbbm{1}_C$. Then $A_i^T\mathbbm{1}_C=1$ for all columns $A_i$ of the matrix $A$ because every perfect matching intersects the tight cut at exactly one edge.
Moreover, $\mathbbm{1}^T \mathbbm{1}_C= |C|$.
Thus, multiplying both sides of $Ax=\mathbbm{1}$ by $\mathbbm{1}_C$ we get $\sum_{i}(A_i^T \mathbbm{1}_C)x_i = \sum_i x_i = |C|$. Notice that the left-hand side of the equation $\sum_i x_i = |C|$ does not depend on the choice of $C$. \hfill $\square$ 
\end{proof}

\begin{definition}
    If a graph $G$ has no non-trivial tight cuts and $G$ is bipartite, then it is a \emph{brace}. If $G$ is not bipartite, then it is a \emph{brick}.
\end{definition}

There is another characterization for bricks and braces, though it will not be relevant in this work:
\begin{theorem}[\cite{lovaszPM}, Section $1$]
    A graph $G$ is a brick if and only if it is $3$-connected and for every pair $x,y$ of vertices, $G-x-y$ has a perfect matching. A bipartite graph $G$ with bipartition $U, W$ with $|U|=|W|$ is a brace if and only if each subset $X \subset U$ with $0<|X|<|U|-1$ has at least $|X|+2$ neighbors in $W$.
\end{theorem}

% The following fact is due to Lov\'asz \cite{lovaszPM}:
% \begin{theorem}
%     The result of contracting one shore of a tight cut of a matching-covered graph is a matching-covered graph.
% \end{theorem}

Existence of tight cuts motivates the following decomposition process on an $r$-graph: \begin{enumerate}
    \item find a tight cut $C$ of $G$, say its shores are $U_1$ and $U_2$,
    \item consider two \textit{$C$-contractions} of $G$: let $G_1=(V\backslash U_1\cup \{u_1\}, E[U_2]\cup C)$ and $G_2=(V\backslash U_2\cup \{u_2\}, E[U_1]\cup C)$, where $E[X]$ denotes the edges of $G$ induced by a subset $X$ of its vertices. In other words, we obtain $G_i$ from $G$ by mapping a shore of $C$ into a single vertex $u_i$ and deleting any loops.
    \item consider each $G_1$ and $G_2$ separately, go back to step $1$.
\end{enumerate}

The process stops when all current graphs in the decomposition have no non-trivial tight cuts i.e. when each graph is either a brick or a brace. This process is called \emph{tight cut decomposition}, and one of its properties is that any implementation of it for a fixed graph $G$ results in the same list of bricks and braces, up to edge multiplicity (\cite{lovaszPM}, Theorem $1.5$).
We now highlight several useful properties of such a decomposition. 
\begin{PR}\label{helpful}
    Let $G=(V,E)$ be a matching covered graph and let $C$ be a tight cut of $G$. Suppose $G_i=(V_i,E_i)$ and $i=1,2$ are two $C$-contractions of $G$. Then \begin{itemize}
        \item[i.] for any perfect matching $M$ of $G$, the set $M\cap E_i$ is a perfect matching of $G_i$, for $i=1,2$;
        \item[ii.] conversely, for any perfect matchings $M_1,M_2$ of $G_1$ and $G_2$, respectively, such that $|M_1\cap M_2|=\{e\}\in C$, the set $M_1\cup M_2$ is a perfect matching of $G$,
        \item[iii.] both $G_1$ and $G_2$ are matching covered;
        \item[iv.] if $G$ is an $r$-graph, then so are $G_1$ and $G_2$.
    \end{itemize} 
\end{PR}
\begin{proof} 
Parts $(i),(ii)$ are immediate from the definition of tight cuts. $(iii)$ follows from $(i)$. 

To prove $(iv)$: any odd cut of $G_i$ is also an odd cut of $G$ and has the same size, thus all odd cuts in $G_i$ have size at least $r$. Lemma \ref{lemma1} guarantees any tight cut must have size exactly $r$, so all vertices in the $C$-contractions of $G$ will have degree $r$.
    \hfill $\square$ 
\end{proof}

In the next sections we discuss how to use this decomposition to inductively construct solutions to $Ax = \mathbbm{1}$ with the desired properties. For convenience, we refer to a (not necessarily nonnegative) vector $x$ satisfying $Ax = \mathbbm{1}$ as a \emph{partitioning} of the edges of $G$ into perfect matchings, and we refer to the matchings corresponding to the support of $x$ as matchings \emph{used} in $x$. Furthermore, we will use interchangeably the notions of a perfect matching and its indicator vector.

\section{Gluing solutions}\label{sect:general_basis_combination}

In this section, we summarize a classical approach (described, e.g., in Murty \cite{Murty}) that combines partitionings of two tight cut contractions of an $r$-graph into a partitioning of that graph.

To this end, let $G=(V,E)$ be an $r$-graph and $C$ be its tight cut. Consider the graphs $G_1$ and $G_2$ obtained from $G$ by contracting one shore of $C$ and let $E_1$, $E_2$ be the edge sets of $G_1$, $G_2$, respectively.
Furthermore, let $A_1$ and $A_2$ be the edge to perfect matching incidence matrices for $G_1$ and $G_2$, respectively.
Given $y,t$ satisfying $A_1y = \mathbbm{1}$ and $A_2t = \mathbbm{1}$, we will construct a vector $x$ satisfying $Ax = \mathbbm{1}$.

Notice that the set of perfect matchings of $G_1$ is a disjoint union of $\{M_i^e\}_{i\in I^e}$, where $e \in C$ and $\{M_i^e\}_{i\in I^e}$ is the set of perfect matchings of $G_1$ containing $e$. Similarly, the set of perfect matchings of $G_2$ is a disjoint union of $\{N_j^e\}_{j\in J^e}$ over all $e\in C$, where $N_j^e$ for $j\in J^e$ are perfect matchings of $G_2$ that contain $e$.
Clearly, all of these matchings do not contain any other edges of $C$.
By part $(ii)$ of Proposition \ref{helpful}, $M_i^e \cup N_j^e =:K_{ij}^e$ is a perfect matching in $G$ for all possible indices $i,j$.
\begin{theorem} \label{thm:gluing}
    For each perfect matching $M$ of $G$, define \begin{equation*}
        x_{M} := \begin{cases}
        y_{M_i^e}t_{N_j^e} & \textmd{ if } M=K^e_{ij}\\
        0 & \textmd{ otherwise}
    \end{cases}.
    \end{equation*} Then $x$ is a partitioning of $G$, i.e. $Ax = \mathbbm{1}_E$.
\end{theorem}
\begin{proof}
     In other words, we need to prove \begin{equation}\label{eq:gluing}
    \sum_{e\in C} \sum_{i \in I^e, j\in J^e}x_{K_{ij}^e}K_{ij}^e=\mathbbm{1}_E.
\end{equation}
It suffices to verify this equality for three types of edges: the ones in $E_1 \backslash C$, in $E_2 \backslash C$, and in $C$.
Every edge $e \in C$ appears on the left side of (\ref{eq:gluing}) with coefficients $y_{M_i^e}t_{N_j^e}$ for all possible $i, j$, so the total sum of the entries corresponding to $e$ in $Ax$ is: \begin{equation*}
    \sum_{i\in I^e,j\in J^e}y_{M_i^e} t_{N_j^e} = (\sum_{i\in I^e} y_{M_i^e})(\sum_{j\in J^e} t_{N_j^e})=1,
\end{equation*} where the last equality follows from $A_1y = \mathbbm{1}_{E_1}$ and $A_2t = \mathbbm{1}_{E_2}$.

Now, projecting the left side of (\ref{eq:gluing}) on the edges in $E_1\backslash C$, note that $K_{ij}^e$ turns into $M_i^e$: \begin{equation*}
    \sum_{e\in C}x_{K^e_{ij}} M_i^e = \sum_{e\in C}( \sum_{i\in I^e}y_{M_i^e}M_i^e\sum_{j\in J^e}t_{N_j^e}) = \sum_{e\in C} \sum_{i \in I^e} y_{M_i^e}M_i^e = \mathbbm{1},
\end{equation*} and similarly for $E_2\backslash C$.

Thus, $x$ corresponds to a partitioning of $G$.  \hfill $\square$ 
\end{proof}

\begin{remark}
    Following the same proof as above, we can show that any assignment of coefficients $x_{K_{ij}^e}$ to $K_{ij}^e$ satisfying the following properties: \begin{equation}\label{eq:partition_required}
        \begin{array}{ccc}
            \sum_{j\in J^e} x_{K_{ij}}^e = y_{M_i^e} & \textmd{ for all } i\in I^e; \\
            \sum_{i\in I^e} x_{K_{ij}}^e = t_{N_j^e} & \textmd{ for all } j\in J^e; \\
            \sum_{i\in I^e, j\in J^e} x_{K_{ij}}^e = 1 & \textmd{ for all } e\in C.
        \end{array}
    \end{equation} will be a partitioning of $G$. This is because the first two conditions ensure that (\ref{eq:gluing}) holds when projected on the edges of $E_1\backslash C$ and $E_2\backslash C$, respectively, and the third one ensures (\ref{eq:gluing}) for the edges of $C$. Furthermore, the third condition is implied by the first two since $\sum_{i\in I^e, j\in J^e} x_{K_{ij}}^e=\sum_{j\in J^e} t_{N_j^e} = 1$ as $t$ corresponds to a partitioning of $G_2$ so the coefficient of $e$ is $1$.
\end{remark}

\section{Main theorem}\label{sect:main}

In this section, we give a different way of combining two solutions for tight cut contractions, which is a crucial ingredient of Theorem \ref{main}. We will construct a solution satisfying the properties (\ref{eq:partition_required}). We will focus on the first two conditions, since the last one is implied. Thus, it suffices to split each $y_{M_i^e}$ into terms $x_{K_{ij}^e}$ while also splitting each $t_{N_j^e}$ into terms $x_{K_{ij}}^e$. A similar approach appears in the proof of the matching polytope theorem in \cite{schrijver}. 

\subsection{A different way to combine solutions}

We begin by proving a weaker version of the main result, which states that there is a partitioning with all coefficients either integral, or equal to $+1/2$. The same algorithm will then be analyzed in more detail in the subsequent section to prove Theorem \ref{main}.

\begin{lemma}\label{lemma2}
There is a solution $x^*$ to $Ax = \mathbbm{1}$ with all entries being either integral, or equal to $+1/2$.
\end{lemma}

\begin{proof}
We prove this by induction on the number of bricks and braces in a tight cut decomposition of the graph. We will prove the base cases (Lemmas \ref{main:brace}, \ref{main:petersen}, and \ref{main:brick}) in the following sections, so for now we assume the statement is true for all bricks and braces.

For the inductive step, consider any graph $G$ with a tight cut $C$ and consider any edge $e \in C$. 
Keep the notations of $A_1, A_2, y, t$ from the previous section. For simplicity, we identify $y_{M_i}^e$ with $y_i$ and $t_{N_j}^e$ with $t_j$ since for the rest of the proof $e$ is fixed. Let further $S_1$ and $S_2$ be the indices at which $y$ and $t$ have negative integers as entries, $P_1$ and $P_2$ be the indices with positive integers, and $H_1$ and $H_2$ be the indices with $+1/2$.
From $A_1 y = \mathbbm{1}$ and $A_2t = \mathbbm{1}$ it follows that \begin{equation*}
    \sum_{i\in P_1}y_i + \sum_{i\in H_1}1/2 - \sum_{i\in S_1}|y_i| = 1 = \sum_{i\in P_2}t_i + \sum_{i\in H_2}1/2 - \sum_{i\in S_2}|t_i|. 
\end{equation*}

Say $-L_1:=\sum_{i\in S_1}y_i$ and $-L_2:=\sum_{i\in S_2}t_i$ (so that both $L_1,L_2$ are nonnegative).
If $L_1\neq L_2$, then without loss of generality $L_1<L_2$ and let us duplicate one matching of $G_1$ from $P_1$ with the largest coefficient $s$ and split it as $s = (s-L_1+L_2) + (L_1-L_2)$, where $s-L_1+L_2$ will now correspond to the nonnegative entries $P_1$ and $L_1-L_2$ will be negative.
This way, the sums of negative coefficients are equal and $\sum_{i\in P_1}y_i + \sum_{i\in H_1}1/2=\sum_{i\in P_2}t_i + \sum_{i\in H_2}1/2$.

Now, we will show how to pair positive and negative coefficients separately.
Suppose we have two sets $a_1, \ldots, a_n$ and $b_1, \ldots, b_m$ (the positive (negative) entries of $y,t$) satisfying $\sum_{i=1}^n a_i = \sum_{j=1}^m b_j=:d$. We may assume both sequences are non-decreasing.
Mark the points of the form $a_1+a_2+\ldots +a_i$ and $b_1+b_2+\ldots+b_j$ (a total of $n+m$ points including $0$ and $d$) on the line segment $[0,d]$.
For every line segment $(w,w')$ of this partition not containing any other partition points, let $i$ be the smallest index with $a_1+\ldots+a_i \geq w'$ and similarly let $j$ be the smallest index for which $b_1+\ldots + b_j \geq w'$. Then let $w'-w$ be the $i+j-1$-th entry of the new vector $c$ (notice that this is well-defined since each subsequent point will have either $i$ or $j$ index increased). 
This $c$ has all entries integral, and its $i+j-1$'st entry will correspond to the union of the perfect matchings associated with $a_i$ and $b_j$. 

We make a small adjustment if there are entries equal to $+1/2$ (by the hypothesis, the only non-integral entries are $+1/2$).
By construction, the fractional coefficients of $a$ and $b$ are listed first. Assuming without loss of generality that $a$ has fewer fractional terms, each $+1/2$ in $a$'s is paired with a $+1/2$ in $b$'s. Until the sequence of $b$'s runs out of $+1/2$'s, all entries of $c$ are $+1/2$ and the remaining entries of $c$ are determined the same way as the all-integral case. Hence, $c$ will have all entries either integral, or equal to $+1/2$.
%This way, the resulting solution (for a fixed $e$) has at most $|H_2|$ coefficients equal to $1/2$, and the rest of the coefficients are integer.
\hfill $\square$
\end{proof}

\subsection{Analysis}

In fact, the above algorithm of combining two solutions preserves several other important properties. Let us prove the following statement. 

\begin{theorem}\label{thm:inductive_step}
    Suppose the method described in Lemma \ref{lemma2} receives as input partitionings $y,t$ of the two $C$ contractions of $G$ and returns a partitioning $x^*$ of $G$. Then, the following properties hold: \begin{enumerate}
        \item[i.] the support size of $x^*$ satisfies $|supp(x^*)|\leq |supp(y)|+|supp(t)|$;
        \item[ii.] the largest entry of $x^*$ satisfies $\|x^*\|_{\infty} \leq \max(\|y\|_\infty,\|t\|_\infty)$;
        \item[iii.] if $y,t$ have all entries in $\mathbb{Z}\cup \{+1/2\}$, then so does $x^*$ and the total number of $+1/2$'s in $x^*$ is at most number of $+1/2$'s in $y$ and $t$ combined;
        \item[iv.] if both $y$ and $t$ only use linearly independent perfect mathcings, then so does $x^*$.
    \end{enumerate}
\end{theorem}
\begin{proof}
    Using the same notation as in the proof of Lemma \ref{lemma2}, for every $e \in C$ the total number of positive terms in $y$ $(t)$ is $H_{1(2)}+P_{1(2)}$ and negative is $S_{1(2)}$. In the proposed algorithm to match $a_i$ and $b_j$, the total number of entries of $c$ is $\leq m+n-1$. Since we might have added an extra term to make $L_1=L_2$, we create a total of at most $\sum_{i=1,2}(H_i+P_i+S_i)$ nonzero entries of $x^*$ corresponding to $e$. Summing this over all $e \in C$ gives $(i)$.    
Similarly, each entry $c_{i+j-1}$ is less than both $a_i$ and $b_j$, thus all entries of $c$ are at most $\min(\max (a_i),\max (b_j))$. Finally, notice that at most one of $a_i$, $b_j$ was artificially augmented by making $L_1=L_2$, so at least one of $\max(a_i),\max(b_j)$ is $\leq \max(\|y\|_\infty, \|t\|_\infty)$, implying $(ii)$.
$(iii)$ holds since for each $e \in C$ the construction yields at most $\max (|H_1|, |H_2|)\leq |H_1|+|H_2|$ terms equal to $+1/2$. Summing this up over all $e$ we get the desired result.    
    
Suppose $(iv)$ does not hold, i.e. there is a linear combination $A w = 0$ where $w$ has the same support as $x^*$. Notice that any column of $A$ is of the form $(M \mid \mathbbm{1}_e \mid N)^T$, where $e=uv \in C$ and $M$, $N$ are some perfect matchings in $G_1-\{u,v\}$ and $G_2-\{u,v\}$.
Comparing projections on $G_1$ and $G_2$, we get that for every matching $M$ of $G_1$ or $G_2$, the sum of entries of $w$ that correspond to a column of $A$ that uses $M$, is zero by assumption on $y,t$.
Hence, it means that for any $e \in C$, the set of perfect matchings containing $e$ has linear combination equal to $0$. Now it suffices to prove that the perfect matchings used in $x^*$ that correspond to the same $e \in C$ are linearly independent.
Indeed, notice that each new matching (going over the matchings in order of the coefficients $c$) either uses a new matching of $G_1$ or of $G_2$, thus does not belong to the linear hull of the previous matchings. 
    %Indeed, for the very first matching we have the interval $(0,c)$ that must correspond to either $a_1$ or $b_1$, so this is the only place where the corresponding matching is used, so the coefficient must be zero.
    %Next, using the induction step, assume that all intervals up to $c_j$ have coefficient zero already.
    %Take the next interval $(c_j, c_{j+1})$. Because $c_{j+1}$ must correspond to either $a$ or $b$ means that there is a matching that is used on this interval and then never used before.
    %Because all of the previous intervals already have coefficient $0$ and the total coefficient for this matching should be zero, it means that the coefficient of $(c_j, c_{j+1})$ must be zero. 
\hfill $\square$
\end{proof}

\subsection{Bipartite case}
To complete the proof of Lemma \ref{lemma2}, we must consider the base cases. We will separately treat braces, Petersen bricks, and non-Petersen bricks.
The braces resulting in the tight cut decomposition of an $r$-graphs must necessarily be bipartite and $r$-regular, as guaranteed by Proposition \ref{helpful}, part $(iv)$. One can check that all regular bipartite graphs satisfy Hall's condition, and hence always contain a perfect matching. Deleting the corresponding edges from the graph again gives a regular bipartite graph. Thus, repeating this process we obtain a union of disjoint perfect matchings that uses all edges (for a more detailed explanation, one can refer to \cite{bondymurty}, Section $5.2$). We remark that these perfect matchings are disjoint, and hence linearly independent. This means that in fact we obtain a $\{0,1\}$ solution with non-zero entries corresponding to \emph{linearly independent} matchings, so it satisfies the conditions of the Theorem \ref{main}. Therefore, we proved the following base case:
\begin{lemma}\label{main:brace}
    Let $G$ be an $r$-graph and a brace, let $A$ be its edge to perfect matching incidence matrix. Then there is a solution $x^*$ to $Ax = \mathbbm{1}$ satisfying the following conditions: \begin{enumerate}
        \item all entries of $x^*$ are in $\{0,1\}$;
        \item $|supp(x^*)|=r$;
        \item all perfect matchings used in $x^*$ are disjoint, and thus linearly independent.
    \end{enumerate}
\end{lemma}

We will now proceed to the case of bricks, which requires some deeper analysis.

\section{Petersen brick solutions}\label{sect:petersen}

In this section, we prove the following:
\begin{lemma}\label{main:petersen}
    Let $G$ be an $r$-graph whose underlying simple graph is a Petersen graph. Let $A$ be its edge to perfect matching incidence matrix. Then there is a solution $x^*$ to $Ax=\mathbbm{1}$ satisfying the following properties: \begin{enumerate}
        \item all entries of $x^*$ are in $\{0,1/2,1\}$;
        \item $|supp(x^*)|\leq m-n+1$ where $m$ and $n$ are the numbers of edges and vertices in $G$, respectively;
        \item $x^*$ has at most $6$ entries equal to $+1/2$;
        \item all perfect matchings used in $x^*$ are linearly independent.
    \end{enumerate}
\end{lemma}
\begin{proof}
    Following the proof of Lemma \ref{lemma2}, we begin with an $r$-graph $G$ and apply tight cut decomposition to it.
Let $B$ be any brick in the tight cut decomposition of $G$ whose underlying graph is the Petersen graph. Notice that $B$ might have parallel edges and it must be $r$-regular.
For convenience, consider a simple graph $B'$, isomorphic to the Petersen graph, whose edges are assigned weights equal to the corresponding number of parallel arcs in $B$.

Recall, Proposition \ref{helpful} guarantees that $B'$ is an $r$-graph. First, let us see which weight assignments for the edges of $B'$ are feasible to ensure this property.
Let $c_{i}$ be the weight of the edge $i$, then we are looking for the $16$-tuples $(c_1, c_2, \ldots, c_{15}, r)$ of positive integers satisfying $deg(v)=r$ for every vertex $v$. This gives $10$ constraints, which can be checked to be linearly independent.
Thus, the set of solutions has dimension $6$ (here, we drop the positive integer requirement for $c_i$'s). 
On the other hand, we can construct a dimension $6$ set of solutions by taking all possible linear combinations of the $6$ perfect matchings of the Petersen graph.
Therefore, any edge weight assignment $c = (c_1,c_2,\ldots, c_{15})$ can be represented as $c = \sum_{i=1}^6 \alpha_i M_i$ where $M_i$ are the perfect matchings of the Petersen graph.

\begin{figure}
    \centering
    \includegraphics[scale=0.4]{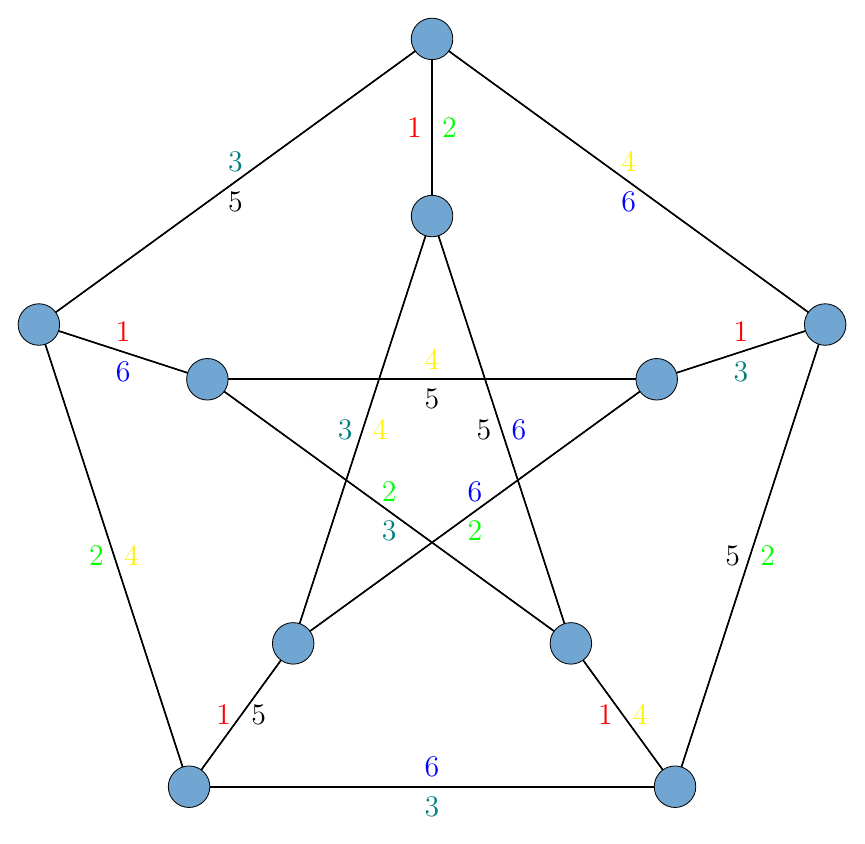}
    \caption{Edges of the Petersen graph labeled by the perfect matchings they belong to.}
    \label{fig:enter-label}
\end{figure}

\begin{claim}
     If $c \in \mathbb{Z}_{> 0}^{15}$, then $\alpha_i \geq 0$ and either $\alpha_i \in \mathbb{Z}$ for all $i$, or $\alpha_i\in 1/2+\mathbb{Z}$ for all $i$.
\end{claim}
\begin{proof}
    Indeed, notice that the total weight of any edge is $\alpha_i+\alpha_j$ for some matchings $M_i$ and $M_j$ because every edge belongs to exactly two perfect matchings. Hence these indices satisfy $\alpha_i \equiv -\alpha_j \mod 1$.
Moreover, any two perfect matchings intersect in one edge, meaning that for any two $i,j$ there is an edge whose weight is exactly $\alpha_i+\alpha_j$, so the previous congruence holds for any pair $i\neq j$ of indices. 
This implies that for any distinct $i,j,k$ we have $\alpha_i \equiv -\alpha_k \equiv \alpha_j$, so $\alpha_i \equiv \alpha_j$.
This together with $\alpha_i \equiv -\alpha_j$ implies $2\alpha_i \equiv 0$, so either $\alpha_i$ is integer and hence all other $\alpha_j$ are integer, or $\alpha_i \equiv 1/2 \mod 1$.

Removing the edges corresponding to $M_1$ from $B$ leaves us with a union of two vertex-disjoint $5$-cycles, meaning that $M_1$ is an odd cut.
Because $B$ is an $r$-graph, every odd cut must have size at least $r$.
The size of this cut is the sum of weights of the edges used in $M_1$, which is $5\alpha_1+\alpha_2+\ldots+\alpha_6$ since every edge in $M_1$ belongs to one other matching.
Hence $r \leq 5\alpha_1+\alpha_2+\ldots+\alpha_6.$
On the other hand, the degree of each vertex is exactly $\sum_i \alpha_i = r$, so \begin{equation*}
    \alpha_1+\alpha_2+\ldots+\alpha_6=r \leq 5\alpha_1+\alpha_2+\ldots+\alpha_6,
\end{equation*} implying $\alpha_1 \geq 0$.
Similarly, all other $\alpha_i$ are non-negative. \hfill $\square$
\end{proof}

This proves that the graph $B'$ has its edge weights $c = \sum_{i=1}^6\alpha_iM_i$ for $\alpha_i\geq 0$ either all integer, or all half-integer.
Now, we show how to return to the original brick $B$ with parallel edges.

\begin{itemize}
    \item If all $\alpha_i$ are integer, write $\sum_{i=1}^6 \alpha_i M_i$ as a sum of $\sum_{i=1}^6 \alpha_i$ single terms. Then, for each edge $e$ of $B'$, its weight $c(e)$ equals the total number of terms using $e$. This means one can replace the instances of $e$ with $e_1, \ldots, e_{c(e)}$, which are the parallel arcs in $B$ corresponding to $e$. Repeating this for every $e \in B'$, we obtain a set of perfect matchings that use each edge of $B$ exactly once, as wanted. Moreover, all perfect matchings will have positive coefficients, so the coefficients are in $\{0, 1\}$. %then for every edge $e$ of $B'$ consider the two matchings using $e$, say they are $M_1$ and $M_2$.
    %Write $\alpha_1M_1+\alpha_2 M_2$ as a sum of $\alpha_1+\alpha_2$ different terms.
    %Notice that $\alpha_1+\alpha_2$ is exactly the weight of $e$ and thus the number of parallel edges in $B$ corresponding to $e$.
   % Denote these edges by $e_1, e_2, \ldots, e_{\alpha_1+\alpha_2}$.
    %For every $i$, replace $e$ with $e_i$ in the $i$-th term of the expanded sum.
   % Clearly, each $e_i$ is now covered exactly once in the new sum.
    %Repeat the same procedure for the remaining edges of $B'$, thus obtaining a set of matchings that covers each edge of $B$ exactly once.
    \item If all $\alpha_i$ are half-integer.
    Let $\beta_i = \alpha_i-1/2 \in \mathbb{Z}_{\geq 0}$.
    Cover an arbitrary subset of edges forming a Petersen graph by $1/2\cdot \sum_{i=1}^6 M_i$, and the remaining is reduced to the previous case for integer $\beta_i$'s. Notice that here, all coefficients must be in $\{0, 1/2, 1\}$ since each edge is covered exactly once and all perfect matchings have nonnegative factors.
\end{itemize}

It remains to prove the bound on support and linear independence.
Notice that both facts are true for Petersen graph with no parallel edges: there are exactly $6$ perfect matchings, which are all linearly independent and their number meets the upper bound of $15-10+1=6$. Now, suppose we have an $r$-graph with weights $c \in \mathbb{Z}^{15}_{+}$ on the edges of Petersen graph. We can view the process of going from $\sum \alpha_i M_i$ to the final partitioning as the following: for some edge $e$ repeatedly replace one copy of a matching $M$ that contains $e$ with a matching $M'$ that uses $e'$, add $e'$ to the graph and reduce the weight of $e$ by one. As a result, the value $m-n+1$ as well as the number of distinct perfect matchings used both increase by one thus maintaining the property \emph{(2)}. Moreover, this step maintains linear independence since the matching $M'$ we add is the unique perfect matching that uses $e'$, and the rest was already linearly independent. Thus, property \emph{(4)} is also preserved.
% Hence, we have a solution with at most $6$ half-integer coefficients for any Petersen brick.
% In fact, there are either no fractional coefficients or there are exactly $6$ of them. We will again remark that the perfect matchings used in the solution are linearly independent. To see this, notice that the six matchings of the Petersen graph $M_i$ are linearly independent and adding parallel arcs keeps the set linearly independent.
\hfill $\square$
\end{proof}

\section{Non-Petersen brick solutions}\label{sect:non-petersen}
The goal of this section is to prove the following:
\begin{lemma}\label{main:brick}
    Let $G$ be an $r$-graph and a non-Petersen brick. Let $A$ be its edge to perfect matching incidence matrix. Then there is a solution $x^*$ to $Ax=\mathbbm{1}$ satisfying the following properties: \begin{enumerate}
        \item all entries of $x^*$ are in $\mathbb{Z}$;
        \item $|supp(x^*)|\leq m-n+1$ where $n$ and $m$ are number of vertices and edges of $G$, respectively;
        \item all perfect matchings used in $x^*$ are linearly independent;
        \item $\|x^*\|_{\infty}\leq 2^{m-n+1}$.
    \end{enumerate}
\end{lemma}
\begin{proof}
 Let $G$ be a non-Petersen brick and an $r$-graph. We will need the following fact dues to Carvalho, Lucchesi, and Murty \cite{dimension}: \begin{theorem}\label{thm:pm_basis}
    For every non-Petersen brick $G$ with $m$ edges and $n$ vertices, the dimension of its matching lattice is $m-n+1$ and it has an integral basis consisting of perfect matching vectors.
\end{theorem}
Here, the perfect matching lattice is the lattice generated by the incidence vectors of perfect matchings of $G$.
Due to Lov\'asz \cite{lovaszPM}, $\mathbbm{1}$ belongs to the matching lattice of $G$, i.e. it can be expressed as an integral combination of perfect matchings. Theorem \ref{thm:pm_basis} guarantees that $\mathbbm{1}$ can be written as a linear combination of at most $m-n+1$ linearly independent perfect matchings.

Using a theory of ear-decompositions of perfect matching-covered graphs, Carvalho, Lucchesi, and Murty \cite{dimension} also provide a construction for the basis of perfect matchings $M_1, M_2, \ldots, M_d$ such that each $M_{i+1}$ contains an edge $e_{i+1}$ not used by $M_1, \ldots, M_i$ for $1\leq i \leq m+n-2$. 
We will use this property in order to prove the following bound: \begin{claim}
    If $x$ is a solution to $\sum_{i=1}^dM_ix_i = \mathbbm{1}$, then $|x_{d-i}|\leq 2^i$.
\end{claim}

\begin{proof}
    By induction on $i$.
In the base case $i=0$, project $\sum_1^d x_i M_i = \mathbbm{1}$ on $e_d$ to get $x_d = 1$ since $e_d$ is only used in the matching $M_d$.
Now, assuming the statement for all previous $i$, projecting $\sum_1^d x_i M_i = \mathbbm{1}$ onto the $e_{d-i-1}$ coordinate, we get: \begin{equation*}
    x_{d-i-1} = 1-(x_{d-i}M_{d-i}+x_{d-i+1}M_{d-i+1}+\ldots + x_d M_d)e_{d-i-1}.
\end{equation*}
Recall that each of $x_{d-j}$ for $j\leq i$ is bounded in absolute value by $2^{j}$, and so are $x_{d-j}M_{d-j}e_{d-i-1}$. Also, $1-x_dM_de_{d-i-1} = 1-M_de_{d-i-1}$ is either $0$ or $1$, so it is bounded in absolute value by $1$. Hence, \begin{equation*}
    |x_{d-i-1}|\leq \sum_{j=1}^i 2^{j} +1 = 2^{i+1},
\end{equation*} as wanted. \hfill $\square$
\end{proof} \hfill $\square$
\end{proof} 

Combining results from before, we get the following statement, which is Theorem \ref{main} together with a bound on the absolute value of the entries of $x^*$: \begin{theorem}
    Let $G$ be an $r$-graph and let $A$ be its edge to perfect matching incidence matrix. Then there is a solution $x^*$ to $Ax = \mathbbm{1}$ satisfying the following conditions: \begin{enumerate}
    \item all non-integral entries of $x^*$ are equal to $+1/2$;
    \item $x^*$ has at most $m-n+1$ non-zero entries;
    \item all non-zero entries of $x^*$ correspond to a linearly independent set of perfect matchings of $G$;
    \item the total number of $+1/2$'s in $x^*$ is at most $6p$, where $p$ is the number of Petersen bricks of $G$;
    \item $||x^*||_{\infty}\leq 2^d$ where $d$ is the largest dimension of the matching lattice (i.e. $m-n+1$) over the non-Petersen bricks of $G$ and $d=0$ if there are no non-Petersen bricks.
\end{enumerate}
\end{theorem}
\begin{proof}
   We prove this by induction on the number of bricks and braces in the tight cut decomposition of $G$. Lemmas \ref{main:brace}, \ref{main:brick}, and \ref{main:petersen} guarantee this for the cases of brick and brace. The inductive step then follows from Theorem \ref{thm:inductive_step}. \hfill $\square$
\end{proof}

\bibliographystyle{splncs04}
\bibliography{PerfectMatching}
%\printbibliography
\end{document}